\documentclass[12pt,twoside,doublespace]{article}

\usepackage{amssymb}
\usepackage{graphicx}

\def\smskip{\par\vskip 5 pt}
\def\QED{\hfill $\Box$\smskip}
\newtheorem{theorem}{Theorem}
\newtheorem{lemma}{Lemma}
\newtheorem{proposition}{Proposition}
\newtheorem{remark}{Remark}

\begin{document}

\begin{center}

\vspace{35pt}

{\Large \bf An Adaptive Partial Linearization Method  }

\vspace{5pt}

{\Large \bf for Optimization Problems on Product Sets}

\vspace{35pt}

{\sc I.V.~Konnov\footnote{\normalsize E-mail: konn-igor@ya.ru}}

\vspace{35pt}

{\em  Department of System Analysis
and Information Technologies, \\ Kazan Federal University, ul.
Kremlevskaya, 18, Kazan 420008, Russia.}

\end{center}

\vspace{35pt}

\begin{abstract}
We suggest an adaptive version of a partial linearization method for
composite optimization problems. The goal function is the sum of
a smooth function and a non necessary smooth convex separable function,
whereas the feasible set is the corresponding
Cartesian product. The method consists in selective component-wise steps
together with a special control of a tolerance sequence.
This technique is destined to reduce the computational expenses per iteration and
maintain the basic convergence properties. We also establish its convergence rates and
describe some examples of applications. Preliminary results of computations
illustrate usefulness of the new method.

\end{abstract}


\newpage


\section{Introduction} \label{s1}

It has been designed a great number of iterative methods for solving various
optimization problems. The custom optimization problem consists in finding
an element in a feasible set $X \subseteq \mathbb{R}^{N}$ that yields
the minimal value of some goal function $\mu : \mathbb{R}^{N} \to \mathbb{R}$ on $X $. For brevity, we write
this problem as
\begin{equation} \label{eq:1.1}
 \min \limits _{\mathbf{x} \in X} \to \mu (\mathbf{x}).
\end{equation}
It is well known that problems with the convex smooth goal function and convex feasible set
constitute one of the most investigated classes in optimization; see e.g. \cite{PD78,Pol83}.
The conditional gradient method is one of the oldest methods in this field. It was first suggested in
\cite{FW56} for the case when the goal function is quadratic and
the feasible set is polyhedral and further was developed by many authors;
see e.g. \cite{LP66,PD78,DR68,Dun80,Pol83}. We recall that the main idea of this method consists in
linearization of the goal function. That is, given the current iterate $\mathbf{x}^{k}\in X$, one finds some
solution $\mathbf{y}^{k}$ of the problem
\begin{equation} \label{eq:1.2}
\min_{\mathbf{y} \in X} \to \langle \mu'(\mathbf{x}^{k}),\mathbf{y}\rangle
\end{equation}
and defines $\mathbf{p}^{k}=\mathbf{y}^{k}-\mathbf{x}^{k}$ as a descent direction at $\mathbf{x}^{k}$.
Taking a suitable stepsize $\lambda_{k} \in (0,1]$, one sets $\mathbf{x}^{k+1}=\mathbf{x}^{k}+\lambda_{k}\mathbf{p}^{k}$
and so on.

During rather long time, this method was not considered as very efficient due to its
relatively slow convergence in comparison with Newton and projection type methods.
However, it has gained a great amount of attention very recently due to several
features significant for many applications, where huge dimensionality and inexact data
create certain drawbacks for more rapid methods. Moreover, in the case of
a polyhedral feasible set its auxiliary problem (\ref{eq:1.2}) appears simpler than those in the
 other methods, and its solution yields usually so-called sparse approximations; see e.g.
\cite{Cla10,Jag13,FG16} and the references therein. It should be noted that
 a great number of applications reduce to problem (\ref{eq:1.1}), where
\begin{equation} \label{eq:1.3}
\mu(\mathbf{x})=f(\mathbf{x})+h(\mathbf{x}),
\end{equation}
$f : \mathbb{R}^{N} \to \mathbb{R}$ is a smooth, but not
necessary convex function, and $h : \mathbb{R}^{N} \to \mathbb{R}$ is
not necessary smooth, but rather simple and convex function.
The appearance of the non-smooth term is caused by
regularization or exact penalty techniques; see e.g.
\cite{Pol83,FSS13}. In this case one can apply the partial linearization (PL for short) method
from \cite{MF81} (see \cite{Pat98,BLM09} for further development), where problem (\ref{eq:1.2}) is replaced with the following:
\begin{equation} \label{eq:1.4}
\min_{\mathbf{y} \in X} \to \langle \mu'(\mathbf{x}^{k}),\mathbf{y}\rangle + h(\mathbf{y}).
\end{equation}
The usefulness of this approach becomes clear if problem
(\ref{eq:1.1}), (\ref{eq:1.3}) is (partially) decomposable, which is typical for very
large dimensional problems.
For instance, let
$$
h(\mathbf{x}) = \sum \limits_{i} h_{i}(\mathbf{x}_{i}) \mbox{ and } \ X = \prod \limits_{i} X_{i}
$$
where $\mathbf{x}_{i} \in X_{i}$. Then (\ref{eq:1.4}) becomes equivalent to several
independent problems of the form
\begin{equation} \label{eq:1.5}
 \min\limits _{\mathbf{y}_{i} \in X_{i}}  \to  \left\{
\left\langle \mathbf{y}_{i}, \frac{\partial f(\mathbf{x}^{k})}{\partial
\mathbf{x}_{i}}\right\rangle +h_{i}(\mathbf{y}_{i})\right\}.
\end{equation}
In case $h \equiv 0$, this decomposition method was considered in \cite{MESS67}.
However, even solution of all the partial problems of form (\ref{eq:1.5})
may appear too expensive. A randomized block-coordinate variant  of the
conditional gradient method was rather recently proposed in \cite{LJSP13}.
A general scheme of block-descent methods for such problems was given in \cite{Pat99}.

We recall for instance that
various engineering problems based on the so-called group LASSO regression method
have this format (see \cite{YL06,MGB08}), as well as many problems of
network resource allocation in wireless multi-user interfering systems
(see \cite{SFSPP14}). We give several additional examples of such decomposable
applied problems in Section \ref{s6}.

The main goal of this paper is to suggest a modification of PL
 methods for decomposable composite optimization problems
of form (\ref{eq:1.1}), (\ref{eq:1.3}), which maintains the basic
convergence properties, but enables one to reduce the computational
expenses per iteration. We follow the approach suggested in
\cite{Kon15d} for regularized splitting methods. The main difference of this method
consists in utilizing PL technique without any regularization in order to simplify the auxiliary problem,
but this implies the dis-continuity of the descent mapping and requires new
substantiation schemes. We take the inexact Armijo type linesearch rule,
which makes our method different from those in \cite{MF81,BLM09}
even in the non-decomposable case.

In what follows, we denote by
$\mathbb{R}^{s}$ the real $s$-dimensional Euclidean space, all
elements of such spaces being column vectors represented by a lower
case Roman alphabet in boldface, e.g. $\mathbf{x}$. We use
superscripts to denote different vectors, and subscripts to denote
different scalars or components of vectors. For any vectors
$\mathbf{x}$ and $\mathbf{y}$ of $\mathbb{R}^{s}$, we denote by $
\langle \mathbf{x}, \mathbf{y} \rangle $ their scalar product, i.e.,
$$
 \langle \mathbf{x}, \mathbf{y} \rangle =\mathbf{x}^{\top} \mathbf{y} =\sum_{i=1}^{s} x_{i}y_{i},
$$
and by $\| \mathbf{x} \|$ the Euclidean norm of $\mathbf{x}$, i.e.,
$\| \mathbf{x}\|=\sqrt{ \langle \mathbf{x}, \mathbf{x} \rangle }$.
We denote by $\mathbb{R}^{s}_{+}$ the non-negative orthant in $\mathbb{R}^{s}$, i.e.
$\mathbb{R}^{s}_{+} = \{ \mathbf{u} \in \mathbb{R}^{s} \ | \ u_{i} \geq 0 \ i=1,\ldots,s \}$.
We also set $\mathcal{R}=\mathbb{R}\bigcup \{-\infty, +\infty\}$.
Given a function $f:\mathbb{R}^{s} \to \mathcal{R}$, we can define its domain
$$
{\rm dom} f=\{ \mathbf{x} \in \mathbb{R}^{s} \ | \  f(\mathbf{x}) > -\infty \}.
$$
For any set $X$, $\Pi (X)$ denotes the family of all nonempty subsets
of $X$.


\section{Problem formulation and preliminary properties}\label{s2}

We first formulate a partitionable optimization problem of form
(\ref{eq:1.1}), (\ref{eq:1.3}). We set
$\mathcal{N}=\{1,\ldots,N\}$ and suppose that there exists a
partition
$$
\mathcal{N}=\bigcup\limits_{i=1}^{n} \mathcal{N}_{i}
$$
with $|\mathcal{N}_{i}|=N_{i}$, $N=\sum \limits_{i=1}^{n} N_{i}$,
and $\mathcal{N}_{i}\bigcap \mathcal{N}_{j}=\varnothing$ if $i\neq
j$ such that
\begin{equation} \label{eq:2.1}
 X=X_{1}\times \dots \times X_{n}=\prod \limits_{i=1}^{n} X_{i},
\end{equation}
where $X_{i}$ is a non-empty, convex, and compact set in
$\mathbb{R}^{N_{i}}$ for $i=1, \dots, n$. Then,  any point
$\mathbf{x}=(x_{1}, \dots,x_{N})^{\top} \in \mathbb{R}^{N}$ is
represented by $\mathbf{x}=(\mathbf{x}_{1}, \dots,
\mathbf{x}_{n})^{\top}$ where $\mathbf{x}_{i}=(x_{j})_{j \in
\mathcal{N}_{i}} \in \mathbb{R}^{N_{i}}$ for $i=1, \dots, n$.
Also, we suppose that
\begin{equation} \label{eq:2.2}
 h(\mathbf{x}) = \sum \limits^{n} _{i=1} h_{i}(\mathbf{x}_{i}),
\end{equation}
where $h_{i}: \mathbb{R}^{N_{i}} \rightarrow \mathcal{R}$ is convex, proper,
lower semi-continuous, and ${\rm dom} h_{i} \supseteq X_{i}$ for $i=1, \dots, n$. Then
the function
$h$ is also convex, proper, and lower semi-continuous and we can define its
subdifferential
$$
\partial h(\mathbf{x})=\partial h_{1}(\mathbf{x}_{1})\times \dots \times \partial
h_{n}(\mathbf{x}_{n}), \quad \forall \mathbf{x}\in X.
$$
So, our problem (\ref{eq:1.1}), (\ref{eq:1.3}),
(\ref{eq:2.1})--(\ref{eq:2.2}) is rewritten as
\begin{equation} \label{eq:2.3}
 \min \limits _{\mathbf{x} \in X_{1}\times \dots \times X_{n}} \to
 \mu(\mathbf{x})=\left\{f(\mathbf{x})+\sum \limits^{n} _{i=1} h_{i}(\mathbf{x}_{i})\right\}.
\end{equation}
Its solution set will be denoted by $X^{*}$ and the optimal value of the
function by $\mu^{*}$, i.e.
$$
\mu^{*} = \inf \limits _{ \mathbf{x} \in X} \mu(\mathbf{x}).
$$
We suppose that the function $f : \mathbb{R}^{N} \to
\mathbb{R}$ is smooth, but not necessary convex. Set
$\mathbf{g}(\mathbf{x})=f'(\mathbf{x})$, then
$$
\mathbf{g}(\mathbf{x})=(\mathbf{g}_{1}(\mathbf{x}), \dots,
\mathbf{g}_{n}(\mathbf{x}))^{\top}, \ \mbox{where} \
\mathbf{g}_{i}(\mathbf{x})=\left(\frac{\partial
f(\mathbf{x})}{\partial x_{j}}\right)_{j \in \mathcal{N}_{i}} \in
\mathbb{R}^{N_{i}}, \ i=1, \dots, n.
$$
From the assumptions above it follows that the function $\mu$ is
directionally differentiable at each point $\mathbf{x} \in X$, that
is, its directional derivative with respect to any vector
$\mathbf{d}$ is defined by the formula:
\begin{equation}\label{eq:2.4}
 \mu '(\mathbf{x}; \mathbf{d})= \langle
{\mathbf{g} (\mathbf{x}), \mathbf{d}}  \rangle + h'(\mathbf{x};
\mathbf{d}), \ \mbox{with} \ h'(\mathbf{x}; \mathbf{d})=\sum
\limits^{n} _{i=1} \max_{\mathbf{b}_{i} \in
\partial h_{i}(\mathbf{x}_{i})} \langle {\mathbf{b}_{i},
\mathbf{d}_{i}}  \rangle ;
\end{equation}
see e.g. \cite{Cla83}.

We need the optimality condition for problem (\ref{eq:2.3}).
\begin{proposition} \label{pro:2.1} \cite[Proposition 2.1]{Kon15d}

(a)  Each solution of problem (\ref{eq:2.3}) is a solution of the
mixed variational inequality (MVI for short): Find a point $\mathbf{x}^{*}
\in X=X_{1}\times \dots \times X_{n}$ such that
\begin{equation} \label{eq:2.5}
\begin{array}{c}
\displaystyle \sum \limits_{i=1}^{n} \left[ \langle  {\mathbf{g}_{i}
(\mathbf{x}^{*}), \mathbf{y}_{i}-\mathbf{x}^{*}_{i}} \rangle   +
h_{i}(\mathbf{y}_{i})-h_{i}(\mathbf{x}^{*}_{i})  \right] \geq 0 \\
\displaystyle \quad \forall \mathbf{y}_{i}\in X_{i}, \quad
\mbox{for} \ i=1, \dots, n.
\end{array}
\end{equation}

(b)  If $f$ is convex, then  each solution of MVI (\ref{eq:2.5})
solves  problem (\ref{eq:2.3}).
\end{proposition}

In what follows, we denote by $X^{0}$ the solution set of MVI
(\ref{eq:2.5})  and call it the set of {\em stationary points} of
problem (\ref{eq:2.3}).

For each point $\mathbf{x}\in X$ we can define a point
$\mathbf{y}(\mathbf{x})=(\mathbf{y}_{1}(\mathbf{x}), \dots,
\mathbf{y}_{n}(\mathbf{x}))^{\top}\in X$ such that
\begin{equation} \label{eq:2.6}
\begin{array}{c}
\displaystyle \sum \limits_{i=1}^{n} \left[ \langle  \mathbf{g}_{i}
(\mathbf{x}),
\mathbf{y}_{i}-\mathbf{y}_{i}(\mathbf{x}) \rangle  +
h_{i}(\mathbf{y}_{i})-h_{i}(\mathbf{y}_{i}(\mathbf{x}))  \right] \geq 0 \\
\displaystyle \quad \forall \mathbf{y}_{i}\in X_{i}, \quad
\mbox{for} \ i=1, \dots, n.
\end{array}
\end{equation}
This MVI gives a necessary and sufficient optimality condition for
the optimization problem:
\begin{equation} \label{eq:2.7}
 \min \limits_{\mathbf{ y}\in X_{1}\times \dots \times X_{n}} \to \sum \limits_{i=1}^{n} \Phi _{i}
(\mathbf{x}, \mathbf{y}_{i}),
\end{equation}
where
\begin{equation} \label{eq:2.8}
\Phi _{i} (\mathbf{x}, \mathbf{y}_{i})= \langle  \mathbf{g}_{i}
(\mathbf{x}),\mathbf{y}_{i}\rangle + h_{i}(\mathbf{y}_{i})
\end{equation}
for $i=1, \dots, n$; cf. (\ref{eq:1.4}). Under the above assumptions
the point $\mathbf{y}(\mathbf{x})$ exists, but is not defined uniquely in general,
hence we can define the set $Y(\mathbf{x})$ of these points at $\mathbf{x}$,
thus defining the set-valued mapping $\mathbf{ x}\mapsto Y(\mathbf{ x})$.
Observe that all the components of $\mathbf{y}(\mathbf{ x})$ can be
found independently, i.e. (\ref{eq:2.7})--(\ref{eq:2.8}) is
equivalent to $n$ independent optimization problems of the form
\begin{equation} \label{eq:2.9}
 \min \limits_{\mathbf{ y}_{i}\in X_{i}} \to \Phi _{i} (\mathbf{x}, \mathbf{y}_{i}),
\end{equation}
for $i=1, \dots, n$ and $\mathbf{y}_{i}(\mathbf{x})$ just solves
(\ref{eq:2.9}). Therefore,
$$
 Y(\mathbf{x}) = Y_{1}(\mathbf{x})\times \dots \times Y_{n}(\mathbf{x}),
$$
where each set $Y_{i}(\mathbf{x})$ is non-empty, convex, and compact.
Moreover, if we set
$$
\sigma _{i} (\mathbf{x}, \mathbf{y}_{i})=\Phi _{i} (\mathbf{x}, \mathbf{x}_{i})-\Phi _{i} (\mathbf{x}, \mathbf{y}_{i})= \langle  \mathbf{g}_{i}
(\mathbf{x}),\mathbf{x}_{i}-\mathbf{y}_{i}\rangle + h_{i}(\mathbf{x}_{i})-h_{i}(\mathbf{y}_{i})
$$
and
$$
\varphi (\mathbf{x})= \sum \limits_{i=1}^{n} \varphi _{i} (\mathbf{x}), \
\varphi _{i} (\mathbf{x})= \max\limits_{\mathbf{ y}_{i}\in X_{i}}\sigma _{i} (\mathbf{x}, \mathbf{y}_{i}) \quad \mbox{for} \ i=1, \dots, n;
$$
then
$$
\varphi _{i} (\mathbf{x})= \sigma _{i} (\mathbf{x}, \mathbf{y}_{i}(\mathbf{x})), \ i=1, \dots, n;
$$
for any $\mathbf{y}(\mathbf{x})=(\mathbf{y}_{1}(\mathbf{x}), \dots,
\mathbf{y}_{n}(\mathbf{x}))^{\top} \in Y(\mathbf{x})$. We can choose the most suitable format for
the definition of a point of $Y(\mathbf{x})$.

We recall that given a set $V \subseteq \mathbb{R}^{s}$,
a set-valued mapping $Q: V \to \Pi (\mathbb{R}^{s})$ is said to be {\em closed} on a set
$W \subseteq V$, if for each pair of sequences $\{ \mathbf{u}^{k} \} \to \mathbf{u}$,   $\{ \mathbf{q}^{k} \} \to \mathbf{q}$
such that $\mathbf{u}^{k} \in W$ and $\mathbf{q}^{k} \in Q(\mathbf{u}^{k})$,
we have $\mathbf{q} \in Q(\mathbf{u})$.

We also need continuity type properties of the marginal functions.


\begin{lemma} \label{lm:2.1}

(a) The function $\varphi : \mathbb{R}^{N} \to
\mathcal{R}$ is lower semi-continuous on $X$;

(b) The mapping $\mathbf{x}\mapsto Y(\mathbf{x})$ is
closed on $X$.
\end{lemma}
{\bf Proof.}
Assertion (a) has been proved in \cite[Lemma 4]{BLM09}. To obtain (b), take sequences
$\{\mathbf{x}^{k}\} \to \bar \mathbf{x}$, $\{\mathbf{y}^{k}\} \to \bar \mathbf{y}$ with
$\mathbf{y}^{k} \in Y(\mathbf{x}^{k})$. Then from (\ref{eq:2.3}) we have
$$
\begin{array}{c}
\displaystyle \sum \limits_{i=1}^{n} \left[ \langle  \mathbf{g}_{i}
(\mathbf{x}^{k}), \mathbf{u}_{i}-\mathbf{y}^{k}_{i} \rangle  +
h_{i}(\mathbf{u}_{i})-h_{i}(\mathbf{y}^{k}_{i})  \right] \geq 0 \\
\displaystyle \quad \forall \mathbf{u}_{i}\in X_{i}, \quad
\mbox{for} \ i=1, \dots, n.
\end{array}
$$
Since $\mathbf{g}$ is continuous and $h$ is lower semi-continuous,
taking the limit $k \to \infty$ gives
$$
\begin{array}{c}
\displaystyle \sum \limits_{i=1}^{n} \left[ \langle  \mathbf{g}_{i}
(\bar \mathbf{x}), \mathbf{u}_{i}-\bar \mathbf{y} \rangle  +
h_{i}(\mathbf{u}_{i})-h_{i}(\bar \mathbf{y})  \right] \geq 0 \\
\displaystyle \quad \forall \mathbf{u}_{i}\in X_{i}, \quad
\mbox{for} \ i=1, \dots, n;
\end{array}
$$
hence $\bar \mathbf{y} \in Y(\bar \mathbf{x})$ and $\mathbf{x}\mapsto Y(\mathbf{x})$
is closed.
\QED

We now show that $\varphi$ can serve as a gap function for problem (\ref{eq:2.3}).


\begin{proposition} \label{pro:2.2}

(a) For any point $\mathbf{x}\in X$ it holds that $\varphi (\mathbf{x})\geq 0$, or, equivalently,
$\varphi_{i} (\mathbf{x})\geq 0$ for $i=1, \dots, n$;

(b) $ \mathbf{x}\in X^{0} \Longleftrightarrow \mathbf{x}\in Y(\mathbf{ x})\Longleftrightarrow
\varphi (\mathbf{x})= 0 \Longleftrightarrow \varphi_{i} (\mathbf{x})= 0, i=1, \dots, n$;
\end{proposition}
{\bf Proof.} Since $\sigma _{i} (\mathbf{x}, \mathbf{x}_{i})=0$, assertion (a) is true.
Next, if $ \mathbf{x}=\mathbf{y}(\mathbf{ x}) \in Y(\mathbf{ x})$, then
(\ref{eq:2.6}) implies $\mathbf{x}\in X^{0}$, $\varphi (\mathbf{x}) \leq 0$ and $\varphi_{i} (\mathbf{x})\leq 0$
for $ i=1, \dots, n$, hence, by (a), $\varphi (\mathbf{x})= 0$ and $\varphi_{i} (\mathbf{x})= 0$
for $ i=1, \dots, n$.  Conversely, let
 $\mathbf{x}$ solve MVI (\ref{eq:2.5}), but $ \mathbf{x}\notin Y(\mathbf{x})$ or $\varphi (\mathbf{x})> 0$.
 Then there exists an index $l$ and a point $\mathbf{x}'_{l} \in X_{l}$ such that
$\sigma _{l} (\mathbf{x}, \mathbf{x}'_{l})>0$. It follows that
$$
-\sum \limits_{i\neq l} \sigma _{i} (\mathbf{x}, \mathbf{x}_{i})- \sigma _{l} (\mathbf{x}, \mathbf{x}'_{l})<0,
$$
i.e. $\mathbf{x} \notin X^{0}$, which is a contradiction. This means that
assertion (b) is true. \QED

We see that the value $\varphi(\mathbf{x})$ can serve as accuracy measure at a point
$\mathbf{x}$. We establish now  a useful descent property. Define for brevity
$I=\{1, \dots, n\}$.


\begin{lemma} \label{lm:2.2} Take
any points $\mathbf{x} \in X$, $\mathbf{y}(\mathbf{x}) \in Y(\mathbf{x})$ and an index $s \in I$. If
$$
\mathbf{d}_{i}= \left\{ {
\begin{array}{ll}
\displaystyle
\mathbf{y}_{s}(\mathbf{x})-\mathbf{x}_{s} \quad & \mbox{if} \ i=s, \\
\mathbf{0} \quad & \mbox{if} \ i\neq s; \\
\end{array}
} \right.
$$
then
\begin{equation} \label{eq:2.11}
 \mu' (\mathbf{x}; \mathbf{d}) \leq -\varphi_{s}(\mathbf{x}).
\end{equation}
\end{lemma}
{\bf Proof.} Due to the definition of $\mathbf{d}$ and (\ref{eq:2.4}), we have
$$
\mu '(\mathbf{x}; \mathbf{d})= \langle {\mathbf{g} (\mathbf{x}),
\mathbf{d}}  \rangle + h'(\mathbf{x}; \mathbf{d})=\langle
\mathbf{g}_{s} (\mathbf{x}), \mathbf{d}_{s}\rangle  +
\max_{\mathbf{b}_{s} \in
\partial h_{s}(\mathbf{x}_{s})} \langle {\mathbf{b}_{s},
\mathbf{d}_{s}}  \rangle.
$$
By convexity, we have
$$
 \langle {\mathbf{b}_{s},
\mathbf{d}_{s}}  \rangle \leq
h_{s}(\mathbf{y}_{s}(\mathbf{x}))-h_{s}(\mathbf{x}_{s})
$$
for any $\mathbf{b}_{s} \in\partial h_{s}(\mathbf{x}_{s})$. It follows that
$$
\mu '(\mathbf{x}; \mathbf{d}) \leq \langle
\mathbf{g}_{s} (\mathbf{x}), \mathbf{y}_{s}(\mathbf{x})-\mathbf{x}_{s} \rangle  +
h_{s}(\mathbf{y}_{s}(\mathbf{x}))-h_{s}(\mathbf{x}_{s})=-\varphi_{s}(\mathbf{x}),
$$
hence (\ref{eq:2.11}) holds true. \QED


\section{The descent method with inexact line-search}\label{s3}

Denote by $\mathbb{Z}_{+}$ the set of non-negative integers. The
basic cycle of the descent PL method with inexact line-search
for MVI (\ref{eq:2.5}) is described as follows.

\medskip
\noindent {\bf Basic cycle (PL).} Choose a point $\mathbf{x}^{0}\in X$
and numbers  $\delta > 0$,  $\beta \in (0,1)$, $\theta
\in (0,1)$.

At the $k$-th iteration, $k=0,1,\ldots$, we have a point
$\mathbf{x}^{k}\in X$.

{\em Step 1:}  Choose an index $s \in I$ such that
$\varphi_{s}(\mathbf{x}^{k}) \geq \delta$, set $s_{k}=s$,
$$
\mathbf{d}^{k}_{i}= \left\{ {
\begin{array}{ll}
\displaystyle
\mathbf{y}_{s}-\mathbf{x}^{k}_{s} \quad & \mbox{if} \ i=s_{k}, \\
\mathbf{0} \quad & \mbox{if} \ i\neq s_{k}; \\
\end{array}
} \right.
$$
where $\mathbf{y}_{s}=\mathbf{ y}_{s}(\mathbf{x}^{k})\in Y_{s}(\mathbf{x}^{k})$ and go to Step 3.
Otherwise (i.e. when $\varphi_{i}(\mathbf{x}^{k}) <
\delta$ for all $i \in I$) go to Step 2.

{\em Step 2:}  Set $\mathbf{z}=\mathbf{x}^{k}$ and stop.

{\em Step 3:}  Determine $m$ as the smallest number in
$\mathbb{Z}_{+}$ such that
\begin{equation} \label{eq:3.1}
 \mu (\mathbf{x}^{k}+\theta ^{m} \mathbf{d}^{k})
 \leq \mu (\mathbf{x}^{k})-\beta \theta ^{m}\varphi_{s}(\mathbf{x}^{k}),
\end{equation}
set $\lambda_{k}=\theta ^{m}$,
$\mathbf{x}^{k+1}=\mathbf{x}^{k}+\lambda_{k}\mathbf{d}^{k}$, and
$k=k+1$. The iteration is complete.
\medskip


\begin{lemma} \label{lm:3.1} The line-search procedure in Step 3 is always finite.
\end{lemma}
{\bf Proof.}
If we suppose that the line-search procedure is infinite, then
$$
\theta^{-m}(\mu (\mathbf{x}^{k}+\theta ^{m} \mathbf{d}^{k}) - \mu
(\mathbf{x}^{k}))>-\beta \varphi_{s}(\mathbf{x}^{k}),
$$
for $m \to \infty$, hence, by taking the limit we have $  \mu'
(\mathbf{x}^{k};\mathbf{d}^{k}) \geq -\beta \varphi_{s}(\mathbf{x}^{k})$,
but Lemma \ref{lm:2.2} gives $ \mu'
(\mathbf{x}^{k};\mathbf{d}^{k}) \leq -\varphi_{s}(\mathbf{x}^{k})$, hence $ (1-\beta
)\varphi_{s}(\mathbf{x}^{k}) \leq  0$, a contradiction. \QED

We recall that a single-valued mapping $\mathbf{p} : \mathbb{R}^{s} \to \mathbb{R}^{s}$
is said to be {\em uniformly continuous} on a set
$V \subset \mathbb{R}^{s}$, if for any number $\varepsilon >0$ there exists a number $\tau >0$
such that $\| \mathbf{p}(\mathbf{x})-\mathbf{p}(\mathbf{y})\| < \varepsilon$ for
each pair of points $\mathbf{x}, \mathbf{y} \in V$ with $\| \mathbf{x}-\mathbf{y}\| < \tau$.
Our convergence analysis will be based on the following property.


\begin{proposition} \label{pro:3.1}
Suppose in addition that the gradient map $\mathbf{g} : \mathbb{R}^{N} \to \mathbb{R}^{N}$
is uniformly continuous on $X$. Then the number of iterations in Basic cycle (PL) is finite.
\end{proposition}
{\bf Proof.}
By construction, we have $- \infty < \mu ^{*}\leq \mu
(\mathbf{x}^{k})$ and $\mu (\mathbf{x}^{k+1})\leq \mu
(\mathbf{x}^{k})-\beta \delta \lambda_{k}$, hence
\begin{equation} \label{eq:3.1a}
\lim \limits_{k\rightarrow \infty }\lambda_{k}=0.
\end{equation}
Besides, the sequence $\{\mathbf{x}^{k}\}$ is bounded and must have limit points,
Suppose that the sequence $\{\mathbf{x}^{k}\}$ is infinite. Since
the set $I$ is finite, there is an index $s_{k}=s$, which is
repeated infinitely. Take the corresponding subsequence $\{k_{l}\}$.
 We intend to evaluate the
difference $\mu (\mathbf{x}^{k_{l}}+\lambda_{k_{l}}\mathbf{d}^{k_{l}}) -
\mu (\mathbf{x}^{k_{l}})$, but we temporarily remove these indices for more convenience.
Then, using the mean value theorem and convexity of $h_{i}$, we have
\begin{eqnarray*}
\displaystyle && \mu (\mathbf{x}+\lambda \mathbf{d})
 - \mu (\mathbf{x}) =f (\mathbf{x}+\lambda \mathbf{d})
 - f(\mathbf{x}) + h_{s}(\mathbf{x}_{s}+\lambda \mathbf{d}_{s})-h_{s}(\mathbf{x}_{s}) \\
  && \leq \lambda \left\{\langle \mathbf{g}_{s}(\mathbf{x}),\mathbf{y}_{s}-\mathbf{x}_{s} \rangle
    + h_{s}(\mathbf{y}_{s})-h_{s}(\mathbf{x}_{s})\right\}
+ \lambda \langle \mathbf{g}_{s}(\mathbf{x}+\xi\lambda \mathbf{d})-\mathbf{g}_{s}(\mathbf{x}),\mathbf{y}_{s}-\mathbf{x}_{s} \rangle \\
 && \leq -\lambda \varphi_{s}(\mathbf{x})+
 \lambda \|\mathbf{g}_{s}(\mathbf{x}+\xi\lambda \mathbf{d})-\mathbf{g}_{s}(\mathbf{x}) \| \|\mathbf{d}_{s} \|,
\end{eqnarray*}
where $\xi=\xi_{k_{l}} \in (0,1)$. Since $X_{s}$ is bounded, $ \|\mathbf{d}_{s} \|\leq C_{s}< \infty$.
Due to the uniform continuity of $\mathbf{g}$, there exists a number $\lambda'>0$
such that
$$
\|\mathbf{g}_{s}(\mathbf{x}+\xi\lambda \mathbf{d})-\mathbf{g}_{s}(\mathbf{x}) \| \leq (1-\beta) \delta/C_{s}
$$
if $\lambda \leq \lambda'$, besides, $\varphi_{s}(\mathbf{x}) \geq \delta$.
It follows that
$$
 \mu (\mathbf{x}+\lambda \mathbf{d})
 - \mu (\mathbf{x}) \leq -\lambda \varphi_{s}(\mathbf{x})+ \lambda (1-\beta) \delta
  \leq -\beta\lambda \varphi_{s}(\mathbf{x})
$$
if $\lambda \leq \lambda'$, hence $\lambda_{k_{l}} \geq \bar \lambda >0$ by the stepsize rule in Basic cycle (PL),
which contradicts (\ref{eq:3.1a}). \QED

The whole method involves the upper level whose iterations (stages)
contain Basic cycle (LP) with decreasing values of $\delta$.

\medskip
\noindent {\bf Method (Upper level).} Choose a point
$\mathbf{z}^{0}\in X$ and a sequence $\{\delta _{l}\} \searrow 0$.

   At the $l$-th stage, $l=1,2,\ldots$, we have a point
$\mathbf{z}^{l-1}\in X$  and a number $\delta _{l}$. Apply Basic
cycle (LP) with $\mathbf{x}^{0}=\mathbf{z}^{l-1}$, $\delta=\delta _{l}$
and obtain a point $\mathbf{z}^{l}= \mathbf{z}$ as its output.
\medskip


\begin{theorem} \label{thm:3.1} Suppose in addition that the gradient
map $\mathbf{g} : \mathbb{R}^{N} \to \mathbb{R}^{N}$
is uniformly continuous on $X$. Then the
sequence $\{\mathbf{z}^{l}\}$ generated by the method with Basic
cycle (LP) has limit points, all these limit points are solutions of
MVI (\ref{eq:2.5}). Besides, if $f$ is convex, then
\begin{equation} \label{eq:3.2}
 \lim \limits_{l\rightarrow \infty} \mu (\mathbf{z}^{l})=\mu^{*};
\end{equation}
and all the limit points of $\{\mathbf{z}^{l}\}$ belong to $X^{*}$.
\end{theorem}
{\bf Proof.}
Following the proof of Proposition \ref{pro:3.1}, we see that
$\mu (\mathbf{z}^{l+1})\leq \mu (\mathbf{z}^{l})$, hence
$$
\lim \limits_{l\rightarrow \infty }\mu (\mathbf{z}^{l})=\tilde \mu.
$$
Besides,  the
sequence $\{\mathbf{z}^{l}\}$ is bounded and must have limit points.
Take an arbitrary limit point $\mathbf{\bar z}$ of
$\{\mathbf{z}^{l}\}$, then
$$
\lim \limits_{t\rightarrow \infty }\mathbf{z}^{l_{t}}=\mathbf{\bar
z}.
$$
For $l>0$ we have
$$
 \varphi_{i}(\mathbf{z}^{l}) \leq \delta_{l} \ \mbox{for all} \ i \in I,
$$
hence $\varphi(\mathbf{z}^{l}) \leq n \delta_{l}$. Due to
Lemma \ref{lm:2.1}, taking the limit $l=l_{t}\rightarrow
\infty$, we obtain $\varphi(\mathbf{\bar z}) \leq0$ and
$\mathbf{\bar z} \in X$. Due to Proposition
\ref{pro:2.2}, this means that $\varphi(\mathbf{\bar z}) =0$
and that the point $\mathbf{\bar z}$
solves MVI (\ref{eq:2.5}). Next, if $f$ is convex, then by
Proposition \ref{pro:2.1} (b), each limit point of
$\{\mathbf{z}^{l}\}$ solves problem (\ref{eq:2.3}). It follows that
$\tilde \mu=\mu^{*}$ and (\ref{eq:3.2}) holds. \QED

In case $h \equiv 0$, the method is a new decomposable version of the
conditional gradient method.
Although the dimensions $N_{i}$ can be arbitrary, we think that the proposed PL method
may have preferences, in particular, over the method from  \cite{Kon15d}, in case when $N_{i}>1$
and all the sets $X_{i}$ are polyhedrons. Also, it may have preferences
 over the usual conditional gradient and partial linearization methods
if the number of subsets $n$ is rather large.


\begin{remark} \label{rm:3.1} The initial boundedness requirement for the feasible set
$X$ was made in Section \ref{s2} only for more simplicity of exposition and
can be replaced with proper coercivity assumptions. In fact, instead of compactness of each set
$X_{i}$ we can require their closedness and  add e.g. the following conditions.

{\bf (C1)} {\em For each $i\in I $ and for each
sequence $\{\mathbf{u}^{l}_{i}\}$ such that $\mathbf{u}^{l}_{i}\in X_{i}$ and 
$\{\|\mathbf{u}^{l}_{i}\|\}\to \infty$ as $l \to \infty$, we have
$\{h_{i}(\mathbf{u}^{l}_{i})/\|\mathbf{u}^{l}_{i}\|\}\to +\infty$.}

{\bf (C2)} {\em For each
sequence $\{\mathbf{u}^{l}\}$ such that $\mathbf{u}^{l}\in X$ 
and $\{\|\mathbf{u}^{l}\|\}\to \infty$ as $l \to \infty$, we have
$\{\mu(\mathbf{u}^{l})\}\to +\infty$.}

Then {\bf (C1)} provides existence of a solution of  auxiliary problem (\ref{eq:2.7})--(\ref{eq:2.8}),
moreover, the sequence $\{\mathbf{d}^{k}\}$ is bounded if so is $\{\mathbf{x}^{k}\}$.
From {\bf (C2)} it follows that $ \mu ^{*} > - \infty$, problems
(\ref{eq:2.3}) and (\ref{eq:2.5}) have solutions, and that
the sequence $\{\mathbf{x}^{k}\}$ is bounded. Therefore,
all the assertions of Section \ref{s3} remain true.

Also, we supposed that ${\rm dom} h_{i} \supseteq X_{i}$ for $i=1, \dots, n$ only for
more simplicity of exposition. Set
$$
D=\prod \limits_{i=1}^{n} ({\rm dom} h_{i} \bigcap X_{i}).
$$
It suffices to assume $D \neq \varnothing$. Then we should only take the initial point
$\mathbf{z}^{0}\in D$.
\end{remark}


\section{Modifications of the linesearch procedure}\label{s4}

Due to Lemma \ref{lm:3.1} the current Armijo rule in (\ref{eq:3.1})
provides its finite implementation e.g. in comparison with
the one-dimensional minimization rule. This version can also be substantiated
under the same assumptions, but we are interested in developing
line-search procedures that are concordant to the partition of the space given
in Section \ref{s2} and do not require calculation of all the components of the
gradient and new point at each iteration. In fact, rule (\ref{eq:3.1})
involves some shift in one component $\mathbf{x}_{s}$,
but utilizes the value of the cost function
at the trial point. That is, we have to calculate the value of $f$ together
with only one component $h_{s}$.

Let us first consider the {\em convex case} where
the function $f$ is convex. Then,  we can replace (\ref{eq:3.1}) with
the following:
\begin{equation} \label{eq:4.1}
\langle  \mathbf{g}_{s} (\mathbf{x}^{k}+\theta^{m}\mathbf{d}^{k}),\mathbf{d}^{k}_{s}\rangle
 + \theta ^{-m}\left\{ h_{s}(\mathbf{x}^{k}_{s}+\theta ^{m}
\mathbf{d}^{k}_{s})-h_{s}(\mathbf{x}^{k}_{s}) \right\}\leq -\beta
\varphi_{s}(\mathbf{x}^{k}).
\end{equation}
Since the
trial point $\mathbf{x}^{k}+\theta ^{m}\mathbf{d}^{k}$ has the shift
from $\mathbf{x}^{k}$ only in $\mathbf{d}^{k}_{s}$, it can be
implemented independently of other variables.  From
(\ref{eq:4.1}) it now follows that
\begin{eqnarray*}
\displaystyle && \mu (\mathbf{x}^{k}+\theta ^{m} \mathbf{d}^{k})
 - \mu (\mathbf{x}^{k}) =f (\mathbf{x}^{k}+\theta ^{m} \mathbf{d}^{k})
 - f(\mathbf{x}^{k}) +
h_{s}(\mathbf{x}^{k}_{s}+\theta ^{m} \mathbf{d}^{k}_{s})-h_{s}(\mathbf{x}^{k}_{s}) \\
 && \leq \theta ^{m} \langle  \mathbf{g}_{s} (\mathbf{x}^{k}+\theta
^{m}\mathbf{d}^{k}),\mathbf{d}^{k}_{i}\rangle +
h_{s}(\mathbf{x}^{k}_{s}+\theta ^{m}
\mathbf{d}^{k}_{i})-h_{s}(\mathbf{x}^{k}_{s})  \leq -\beta
\theta ^{m} \varphi_{s}(\mathbf{x}^{k}),
\end{eqnarray*}
and (\ref{eq:3.1}) holds true. It easy to see that all the assertions of
Section \ref{s3} remain true for this version.

 Moreover, we can utilize even a pre-defined stepsize in the {\em Lipschitz gradient case}.
Let us suppose that partial gradients of the function $f$ are
Lipschitz continuous, i.e.,
$$
\| \mathbf{g}_{i}
(\mathbf{x}+\mathbf{d}^{(i)})-\mathbf{g}_{i}(\mathbf{x})\| \leq
L_{i}\|\mathbf{d}^{(i)}\|=L_{i}\|\mathbf{d}_{i}\|
$$
for any vector $\mathbf{x}$, where
$$
\mathbf{d}^{(i)}_{j}= \left\{ {
\begin{array}{ll}
\displaystyle
\mathbf{d}_{i} \quad & \mbox{if} \ j=i, \\
\mathbf{0} \quad & \mbox{if} \ j\neq i; \\
\end{array}
} \right.
$$
for $i \in I$ and any vector $\mathbf{d}=(\mathbf{d}_{1},
\dots,\mathbf{d}_{n})^{\top} \in \mathbb{R}^{N}$.
Clearly, this property holds if the gradient of $f$ is Lipschitz continuous with
some constant $L >0$, then $L_{i} \leq L$ for each $i \in I$.
It is known that any function $\phi$ having the Lipschitz continuous gradient
satisfies the inequality
$$
\phi (\mathbf{y}) \leq \phi(\mathbf{x})+\langle
\phi'(\mathbf{x}),\mathbf{y}-\mathbf{x} \rangle +0.5L_{\phi}
\|\mathbf{y}-\mathbf{x} \|^{2};
$$
see \cite[Lemma 1.2]{DR68}. Similarly, for any vectors $\mathbf{x}$
and $\mathbf{d}$, we have
$$
f (\mathbf{x}+\mathbf{d}^{(i)}) \leq f(\mathbf{x})+\langle
\mathbf{g}_{i}(\mathbf{x}),\mathbf{d}_{i} \rangle +0.5L_{i}
\|\mathbf{d}_{i} \|^{2} \quad \forall i \in I.
$$
If $\mathbf{d}_{i}= \mathbf{y}_{i}(\mathbf{x})-\mathbf{x}_{i}$, then
we have
\begin{eqnarray*}
\displaystyle && \mu (\mathbf{x}+\lambda \mathbf{d}^{(i)})
 - \mu (\mathbf{x}) =f (\mathbf{x}+\lambda \mathbf{d}^{(i)})
 - f(\mathbf{x}) + h_{i}(\mathbf{x}_{i}+\lambda \mathbf{d}_{i})-h_{i}(\mathbf{x}_{i}) \\
  && \leq \lambda \left\{\langle \mathbf{g}_{i}(\mathbf{x}),\mathbf{d}_{i} \rangle
    + h_{i}(\mathbf{y}_{i}(\mathbf{x}))-h_{i}(\mathbf{x}_{i})\right\}
+0.5L_{i} \lambda^{2}\|\mathbf{d}_{i} \|^{2} \\
 && \leq -\lambda \varphi_{i}(\mathbf{x})+0.5L_{i} \lambda^{2}\|\mathbf{d}_{i} \|^{2}
   \leq -\beta \lambda \varphi_{i}(\mathbf{x}),
\end{eqnarray*}
if
\begin{equation} \label{eq:4.2}
\lambda \leq \bar \lambda_{(i)}(\mathbf{x})=2(1-\beta) \varphi_{i}(\mathbf{x})/(\|\mathbf{d}_{i} \|^{2}L_{i}).
\end{equation}
It follows that (\ref{eq:3.1}) holds with $\lambda_{k} \geq \min\{1, \theta
\bar \lambda_{(s)}(\mathbf{x}^{k})\}>0$.
Moreover, we can simply set $\lambda_{k} = \lambda_{(s)}(\mathbf{x}^{k})>0$,
and all the assertions of Proposition \ref{pro:3.1} and Theorem \ref{thm:3.1}
remain true for this version.  This modification reduces the computational expenses
essentially since calculations of the goal function values are not necessary and  we can calculate
values of the partial gradients $\mathbf{g}_{i}$ and functions $h_{i}$ only for
necessary separate components. Clearly, the adaptive
PL method admits other modifications and extensions, e.g.
selection of a group of indices in $I$ instead of only one
component.

These opportunities make the method very flexible and suitable for
parallel and distributed computations applicable for very
high-dimensional optimization problems; see e.g.
\cite{BT89,Pat99,Jag13,FSS13}.


\section{Convergence rates}\label{s5}

In this section, we give some convergence rates for the adaptive PL method.
We suppose that all the basic assumptions of Section \ref{s2}
hold, but will also utilize some additional conditions.

We first establish the finite termination property
under the following {\em sharp solution condition}, which modifies those
in \cite[Chapter 7, \S 1, Section 3]{Pol83} and
\cite[Section 2.2]{Kon01}.

There exist a number $\tau > 0$ and a point $\bar \mathbf{x} \in X$ such that
$$
\langle \mathbf{g}(\bar \mathbf{x}),  \mathbf{x} - \bar \mathbf{x}
\rangle + h (\mathbf{x}) - h (\bar \mathbf{x}) \geq \tau  \| \mathbf{x} - \bar \mathbf{x} \| \quad \forall \mathbf{x} \in X.
$$


\begin{theorem} \label{thm:5.1}
Let a sequence $\{\mathbf{z}^{l}\}$ be generated by the method with Basic
cycle (LP). Suppose that the function $f$ is convex,
its gradient is Lipschitz continuous with constant $L < \infty$, and that
the sharp solution condition holds. Then
there exists a stage number $t$ such that $X^{*}= Y(\mathbf{z}^{t})$.
\end{theorem}
{\bf Proof.} First we note that the sharp solution condition implies
$\bar \mathbf{x} \in X^{0}$, and, by convexity, $X^{0}=X^{*}$; see Proposition \ref{pro:2.1}.
Next, if there exists some other point $\tilde \mathbf{x} \in X$,
which provides the sharp solution condition, then, again by convexity, we must have
\begin{eqnarray*}
&& \langle \mathbf{g} (\tilde \mathbf{x}),\bar \mathbf{x}-\tilde \mathbf{x} \rangle +
        h(\bar \mathbf{x})-h(\tilde \mathbf{x})
        \leq  \langle \mathbf{g} (\bar \mathbf{x}),\bar \mathbf{x}-\tilde \mathbf{x} \rangle +
        h(\bar \mathbf{x})-h(\tilde \mathbf{x}) \\
&& \leq -\tau  \| \bar \mathbf{x} - \tilde \mathbf{x} \|<0,
\end{eqnarray*}
which is a contradiction. Hence, $X^{*}=\{\bar \mathbf{x} \}$.
From the sharp solution condition for any point $ \mathbf{x} \in X$ we have
\begin{eqnarray*}
&& \langle \mathbf{g} (\mathbf{z}^{l}),\bar \mathbf{x}-\mathbf{x} \rangle +
        h(\bar \mathbf{x})-h(\mathbf{x}) \\
&& =\langle \mathbf{g} (\bar \mathbf{x}),\bar \mathbf{x}-\mathbf{x} \rangle +
        h(\bar \mathbf{x})-h(\mathbf{x})
        + \langle \mathbf{g} (\mathbf{z}^{l})-\mathbf{g} (\bar \mathbf{x}),\bar \mathbf{x}-\mathbf{x} \rangle \\
&&  \leq -\tau  \| \bar \mathbf{x} - \mathbf{x} \|  + L \| \mathbf{z}^{l}-\bar \mathbf{x}\| \| \bar \mathbf{x} - \mathbf{x} \| \\
&&  = -\tau  \| \bar \mathbf{x} - \mathbf{x} \| (1 - L \| \mathbf{z}^{l}-\bar \mathbf{x}\| ).
\end{eqnarray*}
From Theorem \ref{thm:3.1}  we now have $ \{\| \mathbf{z}^{l}-\bar \mathbf{x}\|\}\to 0$
 as $l\to +\infty$. Hence
$$
\langle \mathbf{g} (\mathbf{z}^{l}),\bar \mathbf{x}-\mathbf{x} \rangle +
        h(\bar \mathbf{x})-h(\mathbf{x})<0 \quad \forall \mathbf{x} \in X, \mathbf{x}\neq \bar \mathbf{x},
$$
for $l$ large enough. It follows that there exists a number $t$
such that $ Y(\mathbf{z}^{t})=\{\bar \mathbf{x} \}$. \QED

In the method, each stage contains a finite number of iterations of
the basic cycle. Therefore, it seems suitable
to derive its complexity estimate, which gives the total amount of
work of the method. We now suppose in addition
that the  function $f$ is convex and its partial gradients satisfy Lipschitz continuity
conditions with constants $L_{i}$ for each $i \in I$. Then it was
shown in Section \ref{s4} that we can take the stepsize
\begin{equation} \label{eq:5.1a}
\lambda_{k} = \lambda_{(s)}(\mathbf{x}^{k})
=2(1-\beta) \varphi_{s}(\mathbf{x}^{k})/(\|\mathbf{d}^{k}_{s} \|^{2}L_{s}) \geq
2(1-\beta) \varphi_{s}(\mathbf{x}^{k})/(\rho^{2}L),
\end{equation}
where
$$
L= \max_{s \in I}L_{s}, \ \rho= \max_{s \in I}\rho_{s}, \ \rho_{s}= {\rm diam} X_{s};
$$
see (\ref{eq:4.2}).
We take the value $\Phi(\mathbf{x})=\mu (\mathbf{x})- \mu ^{*}$ as an
accuracy measure for our method. In other words, given a starting
point $ \mathbf{z}^{0}$ and a number $\varepsilon > 0$, we define
the complexity of the method, denoted by $V(\varepsilon )$,  as the
total number of iterations at $l(\varepsilon )$ stages such that
$l(\varepsilon )$ is the maximal number $l$ with
$\Phi(\mathbf{z}^{l}) \geq \varepsilon$, hence,
\begin{equation} \label{eq:5.1}
 V (\varepsilon ) \leq  \sum ^{l(\varepsilon  )} _{l=1}
     V_{l},
\end{equation}
where $V_{l}$ denotes the total number of iterations at  stage
$l$. We proceed to estimate the right-hand side of (\ref{eq:5.1}).
To change $\delta _{l}$, we apply the geometric rate:
\begin{equation} \label{eq:5.2}
 \delta _{l} = \nu ^{l}\delta_{0},  l=0,1,\ldots;
  \quad \nu \in (0,1), \delta_{0}>0.
\end{equation}
By (\ref{eq:3.1}), we have
$$
 \mu (\mathbf{x}^{k+1})
 \leq \mu (\mathbf{x}^{k})-\beta \lambda_{k} \delta _{l},
$$
hence, in view of (\ref{eq:5.1a}), we obtain
\begin{equation} \label{eq:5.3}
 V_{l} \leq \rho^{2}L \Phi(\mathbf{z}^{l-1})/(2\beta(1-\beta)\delta^{2} _{l}).
\end{equation}
Under the above assumptions, for some $ \mathbf{ x}^{*} \in X^{*}$ it holds that
\begin{eqnarray*}
&& \mu  (\mathbf{z}^{l}) - \mu(\mathbf{x}^{*}) = f  (\mathbf{z}^{l}) - f(\mathbf{x}^{*})
 + h  (\mathbf{z}^{l}) - h(\mathbf{x}^{*}) \\
&& \leq \langle \mathbf{g}(\mathbf{z}^{l}),\mathbf{z}^{l}-\mathbf{x}^{*}
\rangle + h(\mathbf{z}^{l})-h(\mathbf{x}^{*})  \\
&& \leq \max\limits_{\mathbf{ y}\in X} \left\{ \langle  \mathbf{g}(\mathbf{z}^{l}),
\mathbf{z}^{l}-\mathbf{y} \rangle  + h(\mathbf{z}^{l})-h(\mathbf{y})  \right\} \\
&& = \varphi (\mathbf{z}^{l}) \leq n \delta_{l}.
\end{eqnarray*}
Using this estimate in (\ref{eq:5.3}) gives
$$
V_{l} \leq \rho^{2}Ln \delta_{l-1}/(2\beta(1-\beta)\delta^{2} _{l}).
$$
From (\ref{eq:5.2}) it follows that
$$
V_{l} \leq \rho^{2}Ln \nu ^{-l}/(2\beta(1-\beta)\delta_{0}\nu)=(C_{1}/\nu)\nu ^{-l}.
$$
Besides, since  $\varepsilon \leq \Phi(\mathbf{z}^{l}) \leq n \delta_{l}=n \delta_{0}\nu ^{l}$, we have
$$
\nu^{-l(\varepsilon  )} \leq n \delta_{0}/\varepsilon.
$$
Combining both the inequalities in  (\ref{eq:5.1}), we obtain
\begin{eqnarray*}
V (\varepsilon ) && \leq C_{1}\nu^{-1}\sum ^{l(\varepsilon  )}
_{l=1} \nu^{-l}  = C_{1} (\nu^{-l(\varepsilon  )}-1)/(1-\nu)  \\
&& \leq C_{1} ((n \delta_{0}/\varepsilon)-1)/(1-\nu).
\end{eqnarray*}
We have obtained the complexity estimate.


\begin{theorem} \label{thm:5.2}
Let a sequence $\{\mathbf{z}^{l}\}$ be generated by the method with Basic
cycle (LP). Suppose that the function $f$ is convex and its partial
gradients satisfy Lipschitz continuity conditions with constants
$L_{i}$ for each $i \in I$.  Then the method has the complexity
estimate
$$
V (\varepsilon ) \leq C_{1} ((n \delta_{0}/\varepsilon)-1)/(1-\nu),
$$
where $C_{1}=\rho^{2}Ln/(2\beta(1-\beta)\delta_{0})$.
\end{theorem}

We observe that the order of the estimates is similar to that in the
usual conditional gradient methods under the same assumptions; see
e.g. \cite{LP66,PD78,Pol83}.


\section{Some examples of applications}\label{s6}

We intend now to give some examples of applied problems
which reduce to decomposable composite optimization problems
of form (\ref{eq:2.3}), where utilization of
the proposed adaptive PL method may give certain preferences.


\subsection{Selective classification problems}

One of the most popular approaches to data classification is
support vector machine techniques; see e.g. \cite{Bur98,Aga15}.
The simplest linear support vector machine problem for data classification
consists in creating an optimal hyperplane separating two convex hulls
of a collection of known
points $\mathbf{x}^{i} \in \mathbb{R}^{m}$, $i=1, \ldots,l$ attributed to previous data observations
with different labels $ y_{i} \in \{-1,+1\}$, $i=1, \ldots,l$,
where $m$ is the number of features. That is, the distance between the hyperplane and each
convex hull should be as long as possible. This separation of the
feature space enables us to classify new data points.
However, this requirement appears too strong for real problems
where the so-called soft margin approach, which minimizes the penalties for
mis-classification, is utilized. This problem can be formulated
as the optimization problem
$$
 \min \limits _{\mathbf{w} \in \mathbb{R}^{n}} \to (1/p)\|\mathbf{w}\|^{p}_{p}
 + C \sum \limits^{l} _{i=1} L( \langle \mathbf{w}, \mathbf{x}^{i} \rangle; y_{i}),
$$
where $L$ is a loss function and $C > 0$ is a penalty parameter. The usual choice is $L(z; y) = \max\{0; 1-yz\}$
 and $p$ is  either 1 or 2. Taking $p=2$, we can rewrite this problem as
$$
 \min \limits _{\mathbf{w}, \mathbf{\xi}} \to 0.5\|\mathbf{w}\|^{2}
 + C \sum \limits^{l} _{i=1} \xi_{i},
$$
subject to
$$
 1-y_{i}\langle \mathbf{w}, \mathbf{x}^{i} \rangle \leq \xi_{i}, \ \xi_{i} \geq 0, \ i=1, \ldots,l.
$$
In this formulation, each observation $i$ is attributed to some data point $\mathbf{x}^{i}$, however,
it seems worthwhile to use sets here since any object may be often represented by some
set of features, this is also the case for noisy observations;
see \cite{SKP08}. So, let object $i$ be represented by a set $X_{i} \in \mathbb{R}^{m}$.
Suppose it is the convex hull of the points $\mathbf{x}^{ik}$, $k=1, \ldots,t$, which thus have the same
label $ y_{i} \in \{-1,+1\}$.
Then we write the soft margin classification problem as follows:
$$
 \min \limits _{\mathbf{w}, \mathbf{\xi}} \to 0.5\|\mathbf{w}\|^{2}
 + C \sum \limits^{l} _{i=1} \xi_{i},
$$
subject to
\begin{eqnarray*}
  1-y_{i}\langle \mathbf{w}, \mathbf{x}^{ik} \rangle \leq \xi_{i},  && \ k=1, \ldots,t; \ i=1, \ldots,l; \\
 \xi_{i} \geq 0, &&  \ i=1, \ldots,l;
\end{eqnarray*}
which somewhat differs from those in \cite{SKP08}.
Here $\xi_{i}$  is the set slack variable
and we impose the penalty for the sum of set slacks. By using the
convex optimization theory, we can write its dual that
has the quadratic programming format:
$$
 \max\limits _{\mathbf{\alpha}} \to \sum \limits^{l} _{i=1} \sum \limits^{t} _{k=1} \alpha_{ik} -
 0.5 \left\| \sum \limits^{l}_{i=1}\sum \limits^{t} _{k=1} (\alpha_{ik}y_{i})\mathbf{x}^{ik}\right\|^{2}
$$
subject to
\begin{eqnarray*}
  && \sum \limits^{t} _{k=1} \alpha_{ik} \leq  C, \ i=1, \ldots,l; \\
 && \alpha_{ik} \geq 0, \ k=1, \ldots,t; \ i=1, \ldots,l;
\end{eqnarray*}
The basic solution of the primal problem is given by the formula:
$$
\mathbf{w}= \sum \limits^{l}_{i=1}\sum \limits^{t} _{k=1} (\alpha_{ik}y_{i})\mathbf{x}^{ik}.
$$
At the same time, we observe that the dual problem falls into
format (\ref{eq:2.3}) and its feasible set is the corresponding
Cartesian product. Hence, our method can be suitable in the high-dimensional case,
where the number of sets is also very large.


\subsection{Network equilibrium problems}

Various network  equilibrium problems represent one of the main tools for
evaluation of flows distribution in traffic and communication networks.
We now describe for instance the path flow formulation
of the network equilibrium problem with elastic demands; see e.g.
\cite{Mag84} and references therein.

The model is determined on an oriented graph, each of its arc being associated with some flow
and some cost (for instance, time of delay), which depends on the values of arc flows.
Usually, the number of nodes and arcs is very large for applied problems.

Let us be given a graph with a finite set of nodes $\mathcal{V}$ and a
set of oriented arcs $\mathcal{A}$ which join the nodes so that any
arc $a=(i, j)$ has the origin $i$ and the destination $j$. Next,
among all the pairs of nodes of the graph we extract a subset of
origin-destination (O/D) pairs $\mathcal{M}$ of the form $m=(i \to
j)$. Besides, each pair $m \in \mathcal{M}$ is associated with a
variable flow demand $v_{m}$ and with  the set of paths
$\mathcal{P}_{m}$ which connect the origin and destination for this
pair. We suppose that each  $v_{m}$ is non-negative with some upper bound $\gamma_{m} \leq +\infty$ for $m
\in \mathcal{M}$. Denote by $\tau_{m}$ the minimal path cost for the pair $m$
and suppose that it depends on the flow demand, i.e.
$\tau_{m}= \tau_{m}(v_{m})$. Also, denote by $u_{p}$ the path flow for the path $p$.
Then the feasible set of flows/demands $W$ can be
defined as follows:
$$
W = \prod_{m \in \mathcal{M}} W_{m},
$$
where
$$
 W_{m}= \left\{ \mathbf{w}_{m}=(\mathbf{u}_{m},v_{m})
  \ \vrule \ \begin{array}{c}
   \sum _{p  \in \mathcal{P}_{m}} u_{p} =v_{m}, u_{p} \geq 0, \
     p \in \mathcal{P}_{m},  \\  0 \leq v_{m} \leq \gamma_{m}, \
  \end{array}
\right\}, \ \forall m \in \mathcal{M},
$$
where $\mathbf{u}_{m}=(u_{p})_{ p \in \mathcal{P}_{m}}$.
 Given a flow vector $\mathbf{u} =   (\mathbf{u}_{m})_{ m \in \mathcal{M}}$,
 one can determine the arc flow
$$
f_{a}=     \sum _{m \in \mathcal{M}} \sum _{p  \in \mathcal{P}_{m}}
\alpha _{p a} u_{p}
$$
 for each arc  $a \in \mathcal{A}$, where
$$
\alpha _{p a}  =\cases {
    1         &\textrm{ if arc $a$ belongs to path $p$,} \cr
    0   &\textrm{  otherwise.} \cr
          }
$$
If the vector $ \mathbf{f} =   (f_{a})_{a \in \mathcal{A}}$ of arc flows is
known, one can determine the arc cost $c_{a}(f_{a})$. We suppose for simplicity that
it depends on the arc flow of just this arc.  Usually, arc costs are monotone
increasing functions of arc flows. Then one can compute costs for each path
$p$:
$$
g_{p}(\mathbf{u}) =  \sum _{a  \in \mathcal{A}}   \alpha _{p a}  c_{a}(f_{a}).
$$
We say that
a feasible flow / demand pair $ (\mathbf{u}^{*}, \mathbf{v}^{*} ) \in W$ is an {\em
equilibrium point} if it satisfies the following conditions:
\begin{equation} \label{eq:6.1}
\forall m \in \mathcal{M}, \ \exists \lambda_{m} \ \mbox{such that}
\ g_{p}(\mathbf{u}^{*}) \left\{
\begin{array}{ll}
\geq \lambda_{m}  \quad & \mbox{if} \quad u^{*}_{p} = 0,  \\
=\lambda_{m}   \quad & \mbox{if} \quad u^{*}_{p}  >0,
\end{array}
\right. \quad   \forall p \in \mathcal{P}_{m};
\end{equation}
and
\begin{equation} \label{eq:6.2}
 \tau_{m}(v_{m}^{*}) \left\{
\begin{array}{ll}
\leq \lambda_{m}  \quad & \mbox{if} \quad v^{*}_{m} = 0,  \\
=\lambda_{m}   \quad & \mbox{if} \quad v^{*}_{m} \in (0,\gamma_{m}); \\
\geq \lambda_{m}
 \quad & \mbox{if} \quad v^{*}_{m} =\gamma_{m};
\end{array}
\right. \quad   \forall p \in \mathcal{P}_{m}.
\end{equation}
However, conditions (\ref{eq:6.1})--(\ref{eq:6.2}) determine
equivalently the following VI: Find a pair $
(\mathbf{u}^{*}, \mathbf{v}^{*}) \in W$ such that
\begin{equation} \label{eq:6.3}
\sum _{m \in \mathcal{M}} \sum _{p  \in \mathcal{P}_{m}}
g_{p}(\mathbf{u}^{*}) (u_{p}- u^{*}_{p})-\sum _{m \in \mathcal{M}}
\tau_{m}(v_{m}^{*}) (v_{m}- v^{*}_{m})  \geq 0       \quad \forall (\mathbf{u}, \mathbf{v})
\in W.
\end{equation}
Furthermore, due to the separability of the functions $c_{a}$ and $\tau_m$, their continuity
implies integrability, i.e., then there exist functions
$$
\eta_a(f_a) = \int \limits_{0}^{f_a} c_a(t) dt \ \forall a \in \mathcal{A}, \
\sigma_m(v_m) = \int \limits_{0}^{v_m} \tau_m(t) dt \  \forall m \in \mathcal{M}.
$$
It follows that VI (\ref{eq:6.3}) also gives an optimality condition of the following optimization problem:
\begin{equation}\label{eq:6.4}
\min \limits_{(\mathbf{u}, \mathbf{v}) \in W} \rightarrow
\left\{ \sum\limits_{a \in \mathcal{A}} \eta_a(f_a) -
 \sum\limits_{m \in \mathcal{M}} \sigma_m(v_m)\right\}.
\end{equation}
Hence, each solution of (\ref{eq:6.4}) is a solution to VI (\ref{eq:6.3}),
the inverse assertion is true if the functions $\eta_a$ and $-\sigma_m$ are convex,
this seems rather natural. However, this problem falls into the basic format
(\ref{eq:2.3}) and the suggested PL method can be applied to this problem.

The basic auxiliary problem consists in finding an element
$(\bar\mathbf{u}_{s}, \bar v_{s})=\mathbf{y}_{s}(\mathbf{x}^{k})\in Y_{s}(\mathbf{x}^{k})$
with $\mathbf{x}^{k}=(\mathbf{u}^{k}, \mathbf{v}^{k})$,
which is now corresponds to a solution of the problem
\begin{equation}\label{eq:6.5}
\min \limits_{(\mathbf{u}_{s}, v_{s}) \in W_{s}} \rightarrow
\left\{ \sum _{p  \in \mathcal{P}_{s}}
g_{p}(\mathbf{u}^{k})u_{p} -
  \sigma_s(v_{s})\right\}
\end{equation}
for some selected pair $s \in \mathcal{M}$. This solution can be found with the simple procedure below,
which is based on optimality conditions (\ref{eq:6.1})--(\ref{eq:6.2}).

First we calculate a shortest path $t \in \mathcal{P}_{s}$ for the pair $s$ with the minimal
cost $\tilde \lambda=g_t(\mathbf{u}^k)$.

Case 1.  If $\tau_{s}(0) \le \tilde \lambda$, then set $\bar v_{s}=0$ and
$\bar u_{p}=0$ for all $p \in \mathcal{P}_{s}$, $\lambda_{s} = \tilde \lambda$.
Otherwise go to Case 2.

Case 2. If $\tau_{s}(\gamma_{s}) \ge \tilde \lambda$,
 set $\lambda_{s} = \tilde \lambda$,
 $\bar v_{s}=\gamma_{s}$, $\bar u_{t}=\gamma_{s}$, and
 $\bar u_{p}=0$ for all $p \in \mathcal{P}_{s}$, $p \neq t$.
 Otherwise go to Case 3.

Case 3. We have $\tau_{s}(\gamma_{s}) < \tilde \lambda < \tau_{s}(0)$. By continuity of
$\tau_{s}$, we find the value $\bar v_{s} \in [0, \gamma_{s}]$  such that
$\tau_{s}(\bar v_{s}) =\tilde \lambda$, set $\lambda_{s} = \tilde \lambda$,
$\bar u_{t}=\bar v_{s}$, and $\bar u_{p}=0$ for all $p \in \mathcal{P}_{s}$, $p \neq t$.

We supposed above that each function $\tau_m$ is continuous, i.e.
that each function $\sigma_m$ is smooth. However, the described procedure for
problem (\ref{eq:6.5})
is extended easily to the case where $-\sigma_m$ is convex
and continuous, then $\tau_m$ can be set-valued.
At the same time, we note that the network equilibrium problem with fixed demands
differs only in somewhat simplified formulation
of problem (\ref{eq:6.4}).
Clearly, the described method remains convergent in these cases and seems in general
simpler and more flexible in comparison with the usual conditional gradient
and projection type methods.


\subsection{Penalty method for decomposable optimization problems}

A great number of optimization problems related to large scale systems
are written as follows:
\begin{equation}\label{eq:6.6}
 \max \ \to \ \sum \limits_{i=1}^{n} \langle \mathbf{c}_{i}, \mathbf{x}_{i} \rangle
\end{equation}
subject to
\begin{eqnarray}
 && \sum \limits_{i=1}^{n} A_{i}\mathbf{x}_{i}  =  \mathbf{b}_{0},  \label{eq:6.7} \\
 &&  \mathbf{x}_{i}  \in X_{i} = \{ \mathbf{y} \in \mathbb{R}^{l_{i} }_{+} \ | \ B_{i}\mathbf{y} \leq \mathbf{b}_{i}
   \}, \ i=1,\ldots,n; \label{eq:6.8}
\end{eqnarray}
for instance, it can be attributed to the total income maximization
in a system containing $n$ subsystems (producers),
who utilize common and particular factors. That is,
producer $i$ chooses an output vector $\mathbf{x}_{i}\in \mathbb{R}^{l_{i}}$,
his/her consumption rates are described by an $m_{0} \times l_{i}$
matrix $A_{i}$  of common factors and by an $m_{i} \times l_{i}$
matrix $B_{i}$  of particular factors, whereas the vector
$\mathbf{c}_{i}$ denotes prices of his/her outputs, the vector $\mathbf{b}_{i} \in
\mathbb{R}^{m_{i}}$ (respectively, $\mathbf{b}_{0}\in \mathbb{R}^{m_{0}}$) denotes inventories of
particular (respectively, common) factors; see e.g. \cite{Las70,Pol83}.
Due to its very large dimensionality, a suitable decomposition approach
can be utilized to reduce the computer memory and calculation expenses.
For instance, the price (Dantzig-Wolfe) decomposition principle
 replaces problem (\ref{eq:6.6})--(\ref{eq:6.8}) with its dual defined with the help of
the Lagrangian including only the term associated with the common constraints in (\ref{eq:6.7}).
However, we can also utilize the penalty approach and replace problem (\ref{eq:6.6})--(\ref{eq:6.8}) with
the sequence of auxiliary problems of the form
\begin{equation}\label{eq:6.9}
 \min \ \to \  0.5 \tau \left\| \sum \limits_{i=1}^{n} A_{i}\mathbf{x}_{i}  - \mathbf{b}_{0} \right\|^{2}
 -\sum \limits_{i=1}^{n} \langle \mathbf{c}_{i}, \mathbf{x}_{i} \rangle
\end{equation}
subject to
\begin{equation}\label{eq:6.10}
\mathbf{x}_{i}  \in X_{i}, \ i=1,\ldots,n;
\end{equation}
where $\tau>0$ is a penalty parameter. Clearly, problem (\ref{eq:6.9})--(\ref{eq:6.10}) also
falls into the basic format (\ref{eq:2.3}) and
application of the suggested PL (conditional gradient) method leads
to some other decomposition method for the initial problem (\ref{eq:6.6})--(\ref{eq:6.8}).
In fact, the partial gradient of the cost function at $\mathbf{x}$ is written as follows
$$
\mathbf{g}_{i}(\mathbf{x})=\tau A^{\top}_{i}\left[ \sum \limits_{j=1}^{n} A_{j}\mathbf{x}_{j}  - \mathbf{b}_{0} \right]
  -\mathbf{c}_{i},
$$
hence
$$
\mathbf{g}_{i}(\mathbf{x}^{k}+\theta \mathbf{d}^{k})=\mathbf{g}_{i}(\mathbf{x}^{k})+ \theta \tau A^{\top}_{i}A_{s}\mathbf{d}^{k}_{s},
$$
and we can make shifts only in the selected component $\mathbf{d}^{k}_{s}$
at each iteration.  Besides, in order to find
$\mathbf{y}_{s}^{k}=\mathbf{y}_{s}(\mathbf{x}^{k})\in Y_{s}(\mathbf{x}^{k})$, we have to solve the separate problem
$$
\min \limits_{\mathbf{y}_{s} \in X_{s}} \rightarrow
\langle \mathbf{g}_{s}(\mathbf{x}^{k}), \mathbf{y}_{s} \rangle.
$$
Combining this method with proper regulation of $\tau$, we obtain a sequence
convergent to a solution of (\ref{eq:6.6})--(\ref{eq:6.8}).


\section{Computational experiments}\label{s7}

In order to compare the performance of the presented method
with the usual non-decomposable version we carried
out preliminary series of computational experiments. For simplicity, we took only
the smooth problems, i.e. set $h \equiv 0$.
Hence, we compared the usual conditional gradient method (CGM) from \cite{PD78} and our method which is treated as
 its adaptive version (ACGM). We took the even partition of $\mathbb{R}^{N}$,
i.e., set $N_{i}=t=N/n$ for $i=1, \dots, n$. Next, each set  $X_{i}$ was chosen to be the standard
simplex in $\mathbb{R}^{t}$, i.e.,
$$
X_{i}=\left\{ \mathbf{u}\in \mathbb{R}^{t}_{+} \ \vrule \ \sum \limits_{i=1}^{t}
u_{i}=1 \right\}.
$$
We took $\Delta_{k}=\varphi(\mathbf{x}^{k})$
as accuracy measure and chose the accuracy $0.1$.
We chose the same starting point $(1/t)e$,
where $\mathbf{e}$ denote the vector of units in $\mathbb{R}^{N}$,
and the rule $\delta_{l+1}=\nu \delta_{l}$  with $\nu = 0.5$ for (ACGM). The methods were
implemented in Delphi with double precision arithmetic.
We report the number of iterations (it) and the total number of calculations (cl) of the
partial gradients $\mathbf{g}_{i}$ for attaining the desired accuracy.

In the first series, we took the convex quadratic cost function.
We chose $\mu(\mathbf{x})=f_{1}(\mathbf{x})$ where
$$
f_{1} (\mathbf{x})= 0.5 \langle P\mathbf{x},\mathbf{x} \rangle -\langle \mathbf{q},\mathbf{x} \rangle,
$$
the elements of the matrix $P$ are defined by
$$
p_{ij}= \left\{ {
\begin{array}{rl}
\displaystyle
\sin (i) \cos (j) \quad & \mbox{if} \ i<j, \\
\sin (j) \cos (i) \quad & \mbox{if} \ i>j, \\
\sum \limits_{s\neq i} | p_{is}| +1 \quad & \mbox{if} \ i=j;
\end{array}
} \right.
$$
and elements of the vector $\mathbf{q}$ are defined by $q_{j}=\sin (j)/j$ for all $i,j$.
The results are given in Table \ref{tbl:1}.
\begin{table}
\caption{The numbers of iterations (it) and partial gradients calculations
(cl)} \label{tbl:1}
\begin{center}
\begin{tabular}{|rr|rr|rr|}
\hline
       &{}           &  (CGM) &{} & (ACGM) &{}  \\
\hline
 $N$ &{} $n$         &  it &{} cl  & it &{} cl      \\
\hline
 10 & 5            &  15 & 75  &  9 & 28  \\
\hline
20 & 5            &  50 & 250  &  108 & 189   \\
\hline
50 & 5             &  143 & 715  &  452 & 676  \\
\hline
100 & 5            &  257 & 1285  &  775 & 1161  \\
\hline
50 & 10           &  228 & 2280   &  632 & 1048  \\
\hline
100 & 10              &  $\Delta_{500}=0.11$  500 & 5000  &  $\Delta_{1500}=0.127$ 1500 & 2515  \\
\hline
80 & 20               &  $\Delta_{500}=0.3$ 500 & 10000  &  766 & 1646  \\
\hline
100 & 20               &  $\Delta_{500}=0.367$ 500 & 10000  &  1328 & 2820  \\
\hline
100 & 25               &  $\Delta_{500}=0.4$ 500 & 12500  &  980 & 2346  \\
\hline
100 & 50               &  $\Delta_{500}=0.76$ 500 & 25000  &  236 & 1036  \\
\hline
\end{tabular}
\end{center}
\end{table}
In the second series, we took the composite convex cost function
$$
\mu(\mathbf{x})=f_{1}(\mathbf{x})+f_{2}(\mathbf{x}),
$$
where $f_{1}$ was defined as above
and
$$
f_{2}(\mathbf{x}) =1/(\langle \mathbf{c},\mathbf{x}\rangle+\tau),
$$
where $c_{i}= 2+\sin(i)$  for $ i=1,\ldots,N$
and $\tau=5$. The results are
given in Table \ref{tbl:2}.
\begin{table}
\caption{The numbers of iterations (it) and partial gradients calculations
(cl)} \label{tbl:2}
\begin{center}
\begin{tabular}{|rr|rr|rr|}
\hline
       &{}           &  (CGM) &{} & (ACGM) &{}  \\
\hline
 $N$ &{} $n$         &  it &{} cl  & it &{} cl      \\
\hline
 10 & 5            &  15 & 75  &  10 & 32  \\
\hline
20 & 5            &  49 & 245  &  113 & 189   \\
\hline
50 & 5             &  139 & 695  &  475 & 666  \\
\hline
100 & 5            &  240 & 1200  &  779 & 1161  \\
\hline
50 & 10           &  231 & 2310   &  620 & 1003  \\
\hline
100 & 10              &  $\Delta_{500}=0.11$  500 & 5000  &  $\Delta_{1500}=0.125$ 1500 & 2515  \\
\hline
80 & 20               &  $\Delta_{500}=0.3$ 500 & 10000  &  766 & 1674  \\
\hline
100 & 20               &  $\Delta_{500}=0.36$ 500 & 10000  &  1329 & 2920  \\
\hline
100 & 25               &  $\Delta_{500}=0.39$ 500 & 12500  &  1011 & 2350  \\
\hline
100 & 50               &  $\Delta_{500}=0.775$ 500 & 25000  &  236 & 1040  \\
\hline
\end{tabular}
\end{center}
\end{table}
In almost all the cases, (ACGM) showed some preference
over (CGM) in the number of partial
gradients calculations.  At the same time, tuning
parameters of (ACGM) needs further
investigations.


\section{Conclusions}

We described a new adaptive component-wise
method for decomposable composite optimization problems involving
non-smooth functions, where the feasible set is the
Cartesian product. The method consists in selective component-wise steps
together with a special control of tolerance sequences. We showed
that this keeps the convergence properties of the usual PL
one together with reduction of the computational
expenses.  We describe several classes of significant applications
for the new method. The preliminary results of computational tests showed rather
satisfactory convergence.


\section*{Acknowledgement}

This work was supported by the RFBR grant, project No. 16-01-00109a
and by grant No. 297689 from Academy of Finland.



\begin{thebibliography}{99}



\bibitem{PD78}
Pshenichnyi,~B.N., Danilin,~Yu.M.:
Numerical Methods in Extremal Problems. MIR, Moscow (1978)

\bibitem{Pol83}
Polyak,~B.T.: Introduction to Optimization. Nauka, Moscow (1983)
[Engl. transl. in Optimization Software, New York (1987)]

\bibitem{FW56}
Frank,~M., Wolfe,~P.: An algorithm for quadratic programming. Nav.
Res. Logist. Quart. 3, 95--110 (1956)

\bibitem{LP66}
Levitin,~E.S., Polyak,~B.T.: Constrained minimization methods.
 USSR Comput. Maths.  Math. Phys. 6, 1--50 (1966)

\bibitem{DR68}
Dem'yanov,~V.F., Rubinov,~A.M.:
Approximate Methods for Solving Extremum Problems. Leningrad University Press, Leningrad (1968)
[Engl. transl. in  Elsevier Science B.V., Amsterdam (1970)]

\bibitem{Dun80}
Dunn, J.C.: Convergence rates for conditional gradient sequences generated by implicit step length
rules. SIAM J. Control Optim. 18, 473--487 (1980)

\bibitem{Cla10}
Clarkson,~K.L.:  Coresets, sparse greedy approximation, and the Frank-Wolfe algorithm.
ACM Trans. on Algor. 6 (4), Art. No. 63, 1--30 (2010)

\bibitem{Jag13}
Jaggi,~M.:  Revisiting Frank-Wolfe: Projection-free sparse convex optimization,  Proc. of the
30th International Conference on Machine Learning (ICML-13), 427--435 (2013)

\bibitem{FG16}
Freund,~R.M., Grigas,~P.:  New analysis and results for the Frank-Wolfe method.
Mathem. Progr. 155,  199--230 (2016)

\bibitem{FSS13}
Facchinei, F., Sagratella, S., Scutari, G.:
Flexible parallel algorithms for big data optimization.  arXiv:1311.2444, 1--8 (November, 11, 2013)

\bibitem{MF81}
Mine,~H., Fukushima,~M.: A minimization method for the sum of a
convex function and a continuously differentiable function. J.
Optim. Theory Appl. 33, 9--23 (1981)

\bibitem{Pat98}
Patriksson,~M.: Cost approximation: a unified framework of descent
algorithms for nonlinear programs. SIAM J. Optim. 8, 561--582 (1998)

\bibitem{BLM09}
Bredies, K., Lorenz, D.A., Maass, P.:  A generalized conditional gradient method and its connection to an
iterative shrinkage method. Comput. Optim. Appl. 42, 173--193 (2009)

\bibitem{MESS67}
Mikhalevich,~V.S., Ermol'ev,~Yu.M., Shkurba,~V.V., Shor,~N.Z.:
Complex systems and the solution of extremal problems.
Kibernetika. 3(5), 29--39 (1967)

\bibitem{LJSP13}
Lacoste-Julien,~S.,  Jaggi,~M., Schmidt,~M., Pletscher,~P.:
Block-coordinate Frank-Wolfe optimization for structural SVMs.
International Conference on Machine Learning, Atlanta (2013) - 31 pp.

\bibitem{Pat99}
Patriksson,~M.:  Nonlinear Programming and Variational Inequality
Problems: A Unified Approach. Kluwer Academic Publishers, Dordrecht
(1999)

\bibitem{YL06}
Yuan, M., Lin, Y.:  Model selection and estimation in regression with
grouped variables. J. R. Statist. Soc. B. 68, 49--67 (2006)

\bibitem{MGB08}
Meier, L., van de Geer, S.,  B\"uhlmann, P.: The group lasso for logistic
regression. J. R. Statist. Soc. B. 70, 53--71 (2008)

\bibitem{SFSPP14}
Scutari,~G.,  Facchinei,~F.,  Song,~P.,  Palomar,~D.P.,  Pang,~J.-S.:
Decomposition by partial linearization: parallel optimization of
multi-agent systems, IEEE Trans. Signal Process.
 62, 641--656 (2014)

\bibitem{Kon15d}
Konnov,~I.V.: Sequential threshold control in descent splitting methods
for decomposable optimization problems. Optim. Meth. Softw. 30, 1238--1254
(2015)

\bibitem{Cla83}
Clarke,~F.H.: Optimization and Nonsmooth Analysis.
John Wiley and Sons, New York (1983)

\bibitem{BT89}
Bertsekas,~D.P., Tsitsiklis,~J.N.: Parallel and Distributed
Computation: Numerical Methods. Prentice-Hall, (1989)

\bibitem{Kon01}
Konnov,~I.V.:   Combined Relaxation Methods for Variational
Inequalities.  Springer, Berlin (2001)

\bibitem{Bur98}
Burges,~C.J.C. A tutorial on support vector machines for pattern recognition.
Data Mining Know. Disc. 2, 121--167 (1998)

\bibitem{Aga15}
Agarval,~C.C.:   Data Mining.  Springer, Heidelberg (2015)

\bibitem{SKP08}
Seref,~O., Kundakcioglu,~O.E., Pardalos,~P.M.: Selective linear and nonlinear classification.
In: Pardalos,~P.M., Hansen,~P. (eds.) Data Mining and Mathematical Programming.
American Mathematical Society, Providence, 211--234 (2008)

\bibitem{Mag84}
Magnanti, T.L.: Models and algorithms for predicting urban traffic equilibria.
In: Florian, M.(ed.) Transportation Planning Models. North--Holland, Amsterdam, 153--185  (1984)

\bibitem{Las70} Lasdon, L.S.:
Optimization Theory for Large Systems.  Macmillan, New York (1970)

\end{thebibliography}
\end{document}